\documentclass[12pt,reqno]{amsart}
\usepackage{amsfonts,amsmath,enumerate}
\usepackage{amssymb}
\usepackage{color}
\numberwithin{equation}{section}
\theoremstyle{plain}
\newtheorem{thm}{Theorem}[section]
\newtheorem{lem}[thm]{Lemma}
\newtheorem{cor}[thm]{Corollary}
\newtheorem{prop}[thm]{Proposition}

\theoremstyle{definition}
\newtheorem{defn}[thm]{Definition}
\newtheorem{ex}[thm]{Example}
\newtheorem{rem}[thm]{Remark}
\newcommand{\al}{\alpha}

\newcommand{\gm}{\gamma}
\newcommand{\Gm}{\Gamma}
\newcommand{\dl}{\delta}

\newcommand{\ep}{\varepsilon}
\newcommand{\et}{\eta}
\newcommand{\tht}{\theta}

\newcommand{\Ld}{\Lambda}
\newcommand{\sg}{\sigma}

\newcommand{\ta}{\tau}

\newcommand{\q}{\quad}

\newcommand{\plimR}{\mathop{\mbox{\rm p-$\lim$}}_{R\to\infty}}
\newcommand{\plim}{\mathop{\mbox{\rm p-$\lim$}}}

\newcommand{\wh}{\widehat}

\renewcommand{\Re}{\mathrm{Re}\,}
\renewcommand{\Im}{\mathrm{Im}\,}

\newcommand{\mA}{\mathcal{A}}
\newcommand{\mC}{\mathcal{C}}
\newcommand{\C}{\mathbb{C}}
\newcommand{\R}{\mathbb{R}}
\newcommand{\N}{\mathbb{N}}

\newcommand{\law}{\mathcal L}

\newcommand{\n}{\noindent}

\newcommand{\bp}{\boxplus}
\newcommand{\xt}{\{X_t, t\ge 0\}}
\newcommand{\yt}{\{Y_t, t\ge 0\}}
\newcommand{\zt}{\{Z_t, t\ge 0\}}

\begin{document}
\setlength{\baselineskip}{18pt}
\setlength{\parindent}{1.8pc}

\title{Selfsimilar free additive processes \\
and freely selfdecomposable distributions}
\renewcommand{\thefootnote}{\fnsymbol{footnote}}
\footnote[0]{Data availability: data sharing not applicable to this article as no datasets were generated.\\
Conflict of interest: to the best of our knowledge, the named authors have no conflict of interest, financial or otherwise.}
\renewcommand{\thefootnote}{\arabic{footnote}}

\author{Makoto Maejima$^1$ and Noriyoshi Sakuma$^2$}

\addtocounter{footnote}{1}\footnotetext{
Keio University, 
E-mail: maejima@math.keio.ac.jp}

\addtocounter{footnote}{1}\footnotetext{
Graduate School of Natural Sciences,
Nagoya City University,
Yamanohata 1, Mizuho-cho, Mizuho-ku, Nagoya, Aichi 467-8501, Japan
E-mail: sakuma@nsc.nagoya-cu.ac.jp\\
}
\date{\today}
\overfullrule=0pt
\maketitle

\vskip 10mm
\begin{abstract}
In the paper by Fan\cite{F06}, he introduced the marginal selfsimilarity of non-commutative stochastic processes and proved the marginal distributions of selfsimilar processes with freely independent increments are freely selfdecomposable.
In this paper, we firstly introduce a new definition, stronger than Fan's one
in general, of selfsimilarity via linear combinations of non-commutative stochastic processes,
although their two definitions are equivalent for non-commutative stochastic processes with freely independent increments.
We secondly prove the converse of Fan's result,
to complete the relationship between selfsimilar free additive processes
and freely selfdecomposable distributions.
Furthermore, we construct stochastic integrals with respect to free additive processes for 
representing the background driving free L{\'e}vy processes of freely selfdecomposable distributions.
A relationship between freely selfdecomposable distributions and their background driving free L{\'e}vy processes in terms of their free cumulant transforms is also given, and several examples are discussed.

\end{abstract}


\vskip 10mm
\section{Introduction}
Free probability theory was introduced by Voiculescu in the 1980s.
He introduced free independence based on the non-commutative probability theory and found applications to operator algebras and random matrices. 
One of the most successful probabilistic objects he introduced is the free additive convolution, denoted by the symbol $\boxplus$, and related analytic tools.
It revealed several mathematical structures
and
enabled us to compute explicit examples in a free probability setting.
We are surprised that there are many interesting correspondences between free and classical probability theory.
Especially, in the 1990s and early 2000, limit theorems and free infinite divisibility were deeply studied and a bijective relationship between free and classical infinitely divisible distributions was found. 
It is, nowadays, called the Bercovici-Pata bijection. 
Infinitely divisible distributions allow us to construct additive processes and L{\'e}vy processes. 
Thus it was natural that many researchers were interested in topics in free probability theory from the point of view of stochastic processes.
Especially, free additive and free L{\'e}vy processes, which were introduced by Biane \cite{B98}, have been studied a lot during the past twenty years.

On the other hand, because of the non-commutativity, it is difficult to consider joint distributions from the view of analytic tools and it is difficult to study non-commutative stochastic processes. 
For example, the selfsimilarity is defined in terms of finite dimensional distributions of stochastic processes in classical probability theory. 
Fan \cite{F06} introduced a weaker notion ``marginal selfsimilarity'' in a free probability setting and gave some free analogue of the results on selfsimilar processes in classical probability theory.
The marginal selfsimilarity is enough whenever we discuss non-commutative 
stochastic processes with freely independent increments, as we will show.
However, for 
further studies of non-commutative stochastic processes not necessarily having
freely independent increments, it might be useful to introduce the selfsimilarity for general non-commutative 
stochastic processes which resemble more closely the classical notion.

The organization of this paper is the following.
In Section 2, we explain free probability, free infinitely divisible distributions, and free additive processes.
In Section 3, a notion of selfsimilarity of non-commutative stochastic processes is introduced, which is stronger than what Fan \cite{F06} defined.
In Section 4, we prove a statement on the relation between selfsimilar free additive processes and freely selfdecomposable distributions, which is a free version of known results in classical probability theory.
In Section 5, we define stochastic integrals with respect to free additive processes, where the Bercovici-Pata bijection plays a crucial role.
In Section 6, we discuss the background driving free L{\'e}vy process of freely selfdecomposable distributions. 
In the final section, we discuss several examples.


\vskip 10mm
\section{Preliminaries: free probability, free infinitely divisible distributions and free additive processes}

In this section, we will gather basic notions and facts, which will be used later.
\subsection{Non-commutative probability space and free independence}

First, we introduce non-commutative probability spaces.
To treat all Borel probability measures, we will use $W^*$-probability space and unbounded operators which are affiliated with such.
For details, see Section $4$ and Appendix in \cite{BT06}.

\begin{defn}
A $W^*$-probability space is a pair $(\mA, \tau)$, where $\mA$ is a von Neumann 
algebra on a Hilbert space $\mathcal H$ and $\tau$ is the faithful, normal, tracial state on $\mA$.
\end{defn}

Let $\mathcal P(\R)$ be the set of all probability measures on $\R$ and let $X$ be a selfadjoint operator in $\mA$.
Then by the Riesz representation theorem, there exists a unique $\mu_X \in \mathcal P (\R)$ such that
$$
\int_{\R} f(x)\mu_X(dx) = \tau(f(X))
$$
for any bounded Borel function $f: \R \to \R$.
We call $X$ a non-commutative random variable.
Throughout this paper, the probability measure of a non-commutative random variable $X$ is written as $\law (X)$.
We only consider selfadjoint operators affiliated with some $W^*$-probability space, so if we say that $X$ is a non-commutative random variable, then it means that it is always selfadjoint.
If $X$ and $Y$ are non-commutative random variables, then $X+Y$ is also a non-commutative random variable.
Selfadjointness is preserved under the addition.

We call that non-commutative random variables $X_{1},X_{2},\ldots,X_{r}$ are freely independent with respect to $\tau$ if for any $k\in\N$, any bounded Borel functions $f_{1}, \ldots, f_{k}:\R\to\R$ and any induces $i_{1},i_{2},\ldots,i_{k}\in\{1,2,\ldots,r\}$ satisfying that $i_{1}\not= i_{2}\not=\dots \not=i_{k}$,
\begin{align*}
\tau\left[ (f_{1}(X_{i_{1}})-\tau(f_{1}(X_{i_{1}}))) 
\cdots  (f_{k}(X_{i_{k}})-\tau(f_{k}(X_{i_{k}})) )
\right] = 0.
\end{align*}
If $X$ and $Y$ are freely independent, $\mu_{X+Y}$ is uniquely determined by $\mu_X$ and $\mu_Y$.
In that case, we call $\mu_{X+Y}$ the free additive convolution of $\mu_X$ and $\mu_Y$ and denote it by $\mu_X\bp\mu_Y$.

In the following, we give analytic tools to compute free additive convolutions.
By $\C^+$ and $\C^-$, we denote the strictly upper and strictly lower complex half planes, respectively.

For $\mu\in\mathcal{P}(\R)$, the Cauchy transform $G_{\mu}: \C^+ \to\C^-$ is defined by
$$
G_{\mu}(z) = \int_{\R} \frac1{z-x}\mu(dx),\q z\in \C^+.
$$
Let $F_{\mu}(z) = \frac{1}{G_{\mu}(z)}: \C^+\to\C^+$.
Then $F_{\mu}(z)$ has the right inverse $F_{\mu}^{-1}(z)$, and the Voiculescu transform $\phi_{\mu}$ is defined by
$$
\phi_{\mu}(z) = F_{\mu}^{-1}(z) - z,\q z\in \Gm (\eta, M),
$$
where $\Gm(\eta, M) = \{ z\in\C \mid |\Re(z)| \le \eta |\Im(z)|, \Im(z)>M\}, \eta>0, M>0.$
Also, the free cumulant transform $\mC_{\mu}$ of $\mu$ is defined by
$$
\mC_{\mu}(z) = z\phi_{\mu}\left(\frac1z\right) = zF_{\mu}^{-1}\left(\frac1z\right)-1,
$$
for $\frac1z \in \Gm(\eta, M)$.

For $\mu_1, \mu_2 \in\mathcal{P}(\R)$, it is known that
$$
\phi_{\mu_1\bp\mu_2} (z) = \phi_{\mu_1}(z) + \phi_{\mu_2}(z),\q 
z\in \Gm (\eta, M),
$$
and
$$
\mC_{\mu_1\bp\mu_2} (z) = \mC_{\mu_1}(z) + \mC_{\mu_2}(z), \q 
z\in \left\{z\,\, \Big| \frac1z \in \Gm(\eta, M)\right\}.
$$
Furthermore, it is known that for $c>0$,
$$
\mC_{D_c\mu}(z) = \mC_{\mu}(cz), \q 
z\in \left\{z\,\, \Big| \frac1z \in \Gm(\eta, M)\right\}.
$$
where $D_c\mu$ is the dilation by $c$, namely $D_c\mu(B) = \mu(c^{-1}B)$ for
all Borel sets $B$ in $\R$.

\subsection{Freely infinitely divisible distributions}
Free additive processes can be constructed by freely infinitely divisible distributions similar to classical additive processes.
In this subsection, we review basic facts on freely infinitely divisible distributions.
\begin{defn}\label{ana-ext}
$\mu\in\mathcal{P}(\R)$ is freely infinitely divisible (also called $\bp$-infinitely divisible), 
if for any $n\in\N$, there exists $\mu_n\in\mathcal{P}(\R)$ such that
$$
\mu = \underbrace{\mu_n\bp\mu_n\bp\cdots\bp\mu_n}_{\mbox{$n$ times}}.
$$
\end{defn}
One necessary and sufficient condition for $\mu\in\mathcal{P}(\R)$ to be freely infinitely divisible is that its Voiculescu transform $\phi_{\mu}(z)$ has an analytic extension defined on $\C^{+}$ with values in $\C^{-}\cup\R$. (See Theorem 5.10 of \cite{BV93}.)
We denote by $ID(\bp)$ the class of all $\bp$-infinitely divisible distributions.
The symbol $ID(*)$ will be used later as the class of all (classical) infinitely divisible distributions.

The following are known (see, e.g. Theorem 4.14 and Proposition 4.16 of \cite{BT06}).

(1)  $\mu\in\mathcal P(\R)$ is $\bp$-infinitely divisible if and only if there exist uniquely a finite measure $\sg$ on $\R$ and $\gm\in\R$ such that
$$
\phi_{\mu}(z) = \gm + \int_{\R} \frac{1+xz}{z-x}\sg(dx),\q z\in\C^+.
$$
The pair $(\gm, \sg)$ is called the free generating pair of $\mu$.

(2) $\mu\in\mathcal P(\R)$ is $\bp$-infinitely divisible if and only if there exist uniquely $a>0, \eta\in\R$ and a free L\'evy measure $\nu$, satisfying $\nu(\{0\})=0$ and $\int_{\R} (|x|^2 \wedge 1)\nu(dx) < \infty$, such that
\begin{equation}\label{cumulant}
\mC _{\mu}(z) = \eta z + az^2 + \int_{\R} \left(\frac1{1-xz} -1-xz{\bf 1}_{[-1,1]}(x)\right ) \nu(dx),\q z\in\C^-.
\end{equation}
The triplet $(a, \nu, \eta)$ is called the free characteristic triplet of $\mu$.

For later use, we recall here the corresponding facts in classical probability theory.

Let $\mu\in ID(*)$.
Then it is known that $\widehat\mu(\theta), \theta\in\R$, the characteristic function of $\mu$,
satisfies
$$
\log \widehat\mu(\theta) = \mathrm{i} \gm\theta +\int_{\R}\left(e^{\mathrm{i} \theta x} -1 - \frac{\mathrm{i} \tht x}{1+x^2}\right) \frac{1+x^2}{x^2}\sg(dx), \q \tht\in\R
$$
and
$$
\log\widehat\mu(\tht) = \mathrm{i} \eta\tht - \frac12a\tht^2 +
\int_{\R}\left( e^{\mathrm{i} \tht x} -1 - \mathrm{i} \tht x {\bf 1}_{[-1,1]}(x)\right)\nu(dx),\q \tht\in\R.
$$
The pair $(\gm, \sg)$ is uniquely determined and is called the generating pair of $\mu$ and the triplet $(a, \nu, \eta)$ is uniquely determined and is called the characteristic triplet of $\mu$.

\subsection{Free additive processes}

\begin{defn}
Let $(\mA , \tau)$ be a $W^*$-probability space and let $\{X_t, t\ge0\}$ be a family of selfadjoint operators affiliated with $\mA$.
We call $\{X_t, t\ge 0\}$ a non-commutative stochastic process on the $W^*$-probability space $(\mA, \tau)$.
\end{defn}

We need the concept of convergence in probability in a free probability setting.

\begin{defn}
Let $X$ and $X_t, t \in [0,\infty)$, be non-commutative random variables.
We say that $X_t\to X$ in probability as $t\to\infty$, if $\law(|X_t-X|)\to\dl_0$ (distribution concentrated at $0$).
We say that $\xt$ is stochastically continuous if for any $t\ge 0$, $X_s \to X_t$ 
in probability as $s\to t$.
\end{defn}

\begin{defn}[Remark 4.28 of \cite{BT06}]
A non-commutative stochastic process $\{X_t, t\ge 0\}$ on some $W^*$-probability space $(\mA, \tau)$ is a free additive process if the following conditions are satisfied.
\begin{enumerate}[\rm (i)]
\item 
For any choice of $n\in\N$ and $0\le t_0 < t_1 <t_2<\cdots <t_n$, the non-commutative random variables
$$
X_{t_0}, X_{t_1}-X_{t_0}, X_{t_2}-X_{t_1}, ... ,X_{t_n}-X_{t_{n-1}}
$$
are freely independent (``freely independent increment property'').
\item 
$X_0=0$.
\item 
$\{X_t,t\ge 0\}$ is stochastically continuous.
\end{enumerate}
\end{defn}

If a free additive process further satisfies that for any $s,t \ge 0$, $\law (X_{s+t}-X_s)$ does not depend on $s$ (``stationary increment property''), a free additive process is called a free L\'evy process.

We can see the following proposition in analogy to classical probability theory.
Note that the marginal distributions of a free L\'evy process are $\bp$-infinitely divisible, (see e.g. p.113 of \cite{BT06}).

\begin{prop}\label{Cconv}
Let $\xt$ be a free L\'evy process.
Then
$$
\mC _{\law(X_t)}(z)= t\mC _{\law(X_1)}(z), \q z \in \C^{-},\, t\ge 0.
$$
\end{prop}

\begin{lem}\label{conv-equiv}
Let $\mu\in ID(\bp)$ and let $\mu_1, \mu_2,...$ be a sequence of probability measures
in $ID(\bp)$ which converge to $\mu\in ID(\bp)$.
Then $\mC_{\mu_n}(z) \to \mC_{\mu}(z)$,  $z \in \C^{-}$, as $n\to \infty$.
\end{lem}

{\it Proof.}
Let $(\gm, \sg)$ and $(a, \nu, \eta)$ be the free generating pair and the free characteristic
triplet of $\mu$, respectively, and
for each $n$, let $(\gm_n, \sg_n)$ and $(a_n, \nu_n, \eta_n)$ be the free generating pair
and the free characteristic triplet of $\mu_n$, respectively.
By Theorem 5.13 of \cite{BT06}, it $\mu_n\to\mu$ as $n\to\infty$, then $\gm_n\to\gm$ 
and $\sg_n\to\sg$ weakly.
On the other hand, if we look at the proof of Proposition 4.16 of \cite{BT06},
we can see that if $\gm_n\to\gm$ and $\sg_n\to\sg$ weakly.
then $a_n\to a, \nu_n\to\nu$ and $\eta_n\to\eta$.
Thus $\mC_{\mu_n}(z)\to \mC_{\mu}(z)$,  $z \in \C^{-}$, as $n\to\infty.$
\qed

\vskip 3mm

{\it Proof of Proposition~\ref{Cconv}.}
In what follows, we use the notation
$$
\bp_{j=1}^n\mu_j = \mu_1\bp \cdots \bp \mu_n \q {\rm for}\q \mu_1,\ldots ,\mu_n\in \mathcal P (\R).
$$
We have, by the stationary increments property, that for $n\ge 1$,
\begin{align*}
\mC _{\law(X_1)}(z) &= \mC_{\law\left(\sum_{j=1}^n(X_{j/n}- X_{(j-1)/n})\right)}(z) \\
& =\mC_{\bp_{j=1}^n \law\left(X_{j/n}- X_{(j-1)/n}\right)}(z) =
n\mC_{\law(X_{1/n})}(z),  \q z \in \C^{-}.
\end{align*}
Hence
$$
\mC_{\law(X_{m/n})} (z)= \frac{m}{n}\mC _{\law(X_1)}(z), \q z \in \C^{-}.
$$
By Lemma~\ref{conv-equiv}, if $\mu_n\in ID(\bp)$ converge to $\mu\in ID(\bp)$ as $n\to\infty$, then $\mC_{\mu_n}(z) \to \mC_{\mu}(z)$, $ z \in \C^{-}$.

If $t>0$ is irrational, choose rational numbers $r_n$ such that $r_n\to t$ as $n\to\infty$.
Since $\{X_t\}$ is stochastically continuous, we have
$$
\mC _{\law(X_{r_n})} (z)\to \mC _{\law(X_t)}(z), \q z \in \C^{-},
$$
as $n\to\infty$, and thus
$$
\mC _{\law(X_t)}(z)= t\mC _{\law(X_1)}(z), \q z \in \C^{-}.
$$
\qed

\vskip 3mm

As mentioned above, the marginal distributions of a free L\'evy process are $\bp$-infinitely divisible.
However, as far as the authors know, the $\bp$-infinite divisibility of marginal distributions of the free additive process is not mentioned anywhere, whose counterpart in classical 
probability theory is known (see e.g., Theorem 9.1 of \cite{S99}).

Hence we start to show it below.

\begin{thm}~\label{stoch-cont}
If $\{X_t, t\ge 0\}$ is a free additive process, affiliated with some $W^*$-probability 
space $(\mA ,\ta )$, then for every $t\ge 0$, $\law(X_t)$ is $\bp$-infinitely divisible.
\end{thm}

\begin{rem}\label{additive-free}
The $\bp$-infinitely divisibility of the marginal distributions of free additive processes is very important for the use of the Bercovici-Pata bijection between $ID(*)$
and $ID(\bp)$, which will be defined in Section 5.
\end{rem}

The proof of Theorem~\ref{stoch-cont} can be carried out in a
similar way to that of the corresponding result in classical probability theory, 
(see again, e.g., Theorem 9.1 of \cite{S99}).

We need the following free-version of Khinchine's theorem in \cite{BP99}.

\begin{prop}\label{BP}
Let $\{\mu_{nj} \,\,| \,\,n\in\N, 1\le j \le k_n\}$ be an array of probability measures on $\R$ and 
$\{c_n, n\in\N\}$ a sequence of real numbers such that for every $\ep>0$
$$
\lim_{n\to \infty}\max_{1\le j\le k_n} \mu_{nj}\{(x : |x|>\ep)\}=0 
$$ 
and the probability measures 
$\mu_n = \dl _{c_n}\bp\mu_{n1}\bp\mu_{n2}\bp\cdots \bp \mu_{nk_n}$
converge weakly to a probability measure $\mu$ on $\R$.
Then $\mu$ is $\bp$-infinitely divisible.
\end{prop}

\begin{lem}
A stochastically continuous non-commutative stochastic process $\{X_t, t\ge 0\}$ is uniformly 
stochastically continuous on any finite interval $[0, T]$, that is, for any $\ep >0$ and $\et>0$,
there is $\dl >0$ such that, if $s$ and $t$ are in $[0, T]$ and satisfy $|s-t|<\dl$, then
$$
\mu_{s,t} (\{ x : |x|>\ep\})<\et,
$$
where $\mu_{s,t} = \law(X_t-X_s)$.
\end{lem}

The proof of this lemma can be carried out in exactly the same way as for that of the
corresponding result in classical probability theory (see Lemma 9.6 of \cite{S99}).

\vskip 3mm

{\it Proof of} Theorem~\ref{stoch-cont}.
Fix $t>0$ and let $t_{n,j}=jt/n$ for $n=1,2,...$ and $j=0,1, ... , n$.
Let $k_n=n$ and $Y_{nj}=X_{t_{n,j}}-X_{t_{n,j-1}}$ for $j=1, ... ,n$.
Then $\{\mu_{nj}=\law(Y_{nj})  | n\ge 1, 1\le j \le n \}$ is an array of probability measures
satisfying all conditions in Proposition~\ref{BP}.
The row sum $\sum_{j=1}^{n}Y_{nj}$ equals $X_t$.
Thus we can apply Proposition~\ref{BP} with $\mu=\law(X_t)$ and $c_n=0$ to conclude the theorem.
\qed

\vskip 3mm

Furthermore, we have the following, which is also important as well as the
$\bp$-infinite divisibility of $\law(X_t)$.
(This is a free version of a part of Theorem 9.7 (i) of \cite{S99} in classical probability theory.)

\begin{thm}\label{dif-id}
If $\{X_t, t\ge 0\}$ is a free additive process, affiliated with some $W^*$-probability 
space $(\mA ,\ta )$, then for every $0\le s\le t<\infty$, $\law(X_t-X_s)$ is $\bp$-infinitely divisible.
\end{thm}

{\it Proof.}
Fix $u\ge 0$.
Then $\{X_{u+t} - X_u, t\ge 0\}$ is a free additive process.
Hence $\law(X_t-X_s)$ is $\bp$-infinitely divisible by Theorem~\ref{stoch-cont}.
\qed


\vskip 10mm
\section{Selfsimilar stochastic processes}

Let us consider selfsimilar processes in a free probability setting.  
In classical probability theory, a selfsimilar process $\xt$ is defined as follows:  
For any $a>0$, there exists $b=b(a)$ such that, for any choice of $n\in\N$ and $0\le t_1<t_2< \cdots <t_n$,
\begin{equation}\label{ss}
\law\left(\left(X_{at_1}, X_{at_2},\ldots , X_{at_n}\right)\right)
= \law\left(b\left(X_{t_1}, X_{t_2},\ldots ,X_{t_n}\right)\right).
\end{equation}
Argumentation based on the finite dimensional distributions of the process cannot be used in free probability theory.
However, in classical probability theory, we have the following.

\begin{prop}\label{SS}
Let $\xt$ be a stochastic process.
Suppose that all linear combinations
$$
\sum_{j=1}^n c_j X_{t_j,}\q n\in\N,\,\, c_1,c_2, ..., c_n\in\R,\,\, 0\le t_1<t_2<\cdots <t_n,
$$
are selfsimilar, namely, for any $a>0$, there exists $b=b(a)$, 
independent of the choice of $\{c_j\}$ and $\{t_j\}$, such that
$$
\law\left( \sum_{j=1}^n c_{j}X_{at_j}\right) = \law\left( b \sum_{j=1}^n c_jX_{t_j}\right).
$$
Then $\xt$ satisfies \eqref{ss}.
\end{prop}

{\it Proof.}
Let 
$$
{\bf X} = (X_{at_1}, X_{at_2}, ... , X_{at_n}) \q {\rm and} \q {\bf Y} =( b X_{t_1}, bX_{t_2}, ... ,b X_{t_n}).
$$
We want to show $\law ({\bf X}) = \law({\bf Y})$.
Since for any ${\bf c}=(c_1, c_2, ... , c_n)$, $\law (({\bf c}, {\bf X})) = \law (({\bf c}, {\bf Y}))$,
we have $E\left[ e^{{\rm i}({\bf c}, {\bf X})}\right] = E\left[ e^{{\rm i}({\bf c}, {\bf Y})}\right]$.
Thus $\bf X$ and $\bf Y$ have the identical characteristic function, implying that
$\law({\bf X}) = \law ({\bf Y})$, namely, \eqref{ss} holds.
Thus $\xt$ is selfsimilar.
\qed

In view of Proposition~\ref{SS}, the following definition would be acceptable in free probability theory.

\begin{defn}\label{free-ss}
Let $(\mA, \tau)$ be a $W^*$-probability space and let $\{X_t, t\ge0\}$ be a family of selfadjoint operators affiliated with $\mA$.
We call $\{X_t,t\ge 0\}$ a selfsimilar non-commutative stochastic process, if for any $a>0$, there exists $b=b(a)$ dependent only on $a$ such that all linear combinations
$$
\law\left( \sum_{j=1}^n c_{j} X_{at_j}\right) = \law\left( b \sum_{j=1}^n c_jX_{t_j}\right),\q
n\in\N, \,\, c_1, c_2, ... , c_n\in \R, \,\, 0\le t_1<t_2< \cdots <t_n.
$$
\end{defn}

This definition defines a smaller class than what Fan \cite{F06} defined.
Namely,  Fan defined the selfsimilar process $\xt$ in the way that
if for any $a>0$, there exists $b>0$ such that
\begin{equation}\label{mss}
\law (X_{at}) = \law (bX_t),\q t\ge 0.
\end{equation}
To distinguish Fan's definition from ours, we call the former for marginal selfsimilarity.
Let $\xt$ be a selfsimilar non-commutative stochastic process in the sense of
Definition~\ref{free-ss}, satisfying $\lim_{t\to 0}X_t = X_0$ in probability. 
Then there exists a unique $H>0$ such that $b(a)=a^H$ for all $a>0$.
This fact is also true for marginal selfsimilar processes as shown in Theorem 3.4 of \cite{F06}, where Fan follows the argument in the proof of Theorem 1.1.1 of \cite{EM02}
in classical probability theory.
This $H$ is called the index of selfsimilarity, and $\xt$ is called $H$-selfsimilar. 

It is obvious that if a non-commutative stochastic process $\xt$ is a selfsimilar process, then $\xt$ is marginally selfsimilar.
However, the converse is not necessarily true.
Actually, in 2 in Remarks at page 312 of \cite{ST94} in classical probability theory, there is an example which is marginally selfsimilar but not selfsimilar.
The same example also works in free probability theory.
This example is that $\law(X_t) = \law(t^H X_1)$ for each $t>0$ and 
it is not freely additive.

However, as far as we consider free additive processes, the marginal selfsimilarity
implies the selfsimilarity in the sense of Definition~\ref{free-ss},
namely, both definitions are the same.
In order to prove this, we first show the following.

\begin{thm}\label{same}
Let $\{X_t, t\ge 0\}$ and $\{X'_t, t\ge 0\}$ be two free additive processes such that
$\law(X_t)=\law(X_t')$ for all $t\ge 0$.
Then 
for any $n\in\N , c_1, c_2, ... ,c_n,\in\R,  0\le t_1 <t_2< \cdots <t_n$,
\begin{eqnarray}\label{lineq}
\law\left( \sum_{j=1}^n c_j X_j\right) = \law\left( \sum_{j=1}^n c_j X'_j\right).
\end{eqnarray}
\end{thm}

{\it Proof.}
Let 
$\mu_{s,t}=\law(X_t-X_s)$ and $\mu'_{s,t}=\law(X'_t-X'_s)$ for $0\le s\le t<\infty$.
Then $\mu_{0,t} = \mu'_{0,t}$ since $X_0 = X_0'$.
By the freely independent increment property, we have
$$
\mC _{\mu_{0,s}}(z) + \mC_{\mu_{s,t}}(z)= \mC _{\mu_{0,t}}(z), \q z \in \C^{-},
$$
and
$$
\mC _{\mu'_{0,s}}(z) + \mC_{\mu'_{s,t}}(z)= \mC _{\mu'_{0,t}}(z), \q z \in \C^{-}.
$$
Thus  $\mC_{\mu_{s,t}}(z) = \mC_{\mu'_{s,t}}(z),$ namely
$$
\mC_{\law(X_t-X_s)}(z) = \mC_{\law(X'_t-X'_s)}(z), \q z \in \C^{-}.
$$
First, we introduce a simple but important equation. 
For any $n\in\N, c_1, c_2, ... ,c_n,\in\R$, 
$0\le t_1 <t_2< \cdots <t_n$,
$$
\sum_{j=1}^n c_j X_{t_j} =
\sum_{k=1}^n\left(\sum_{j=n+1-k}^n c_j\right)\left(X_{t_{n+1-k}}-X_{t_{n-k}}\right),\q (t_0=0).
$$
Put $c_{n,k} = \sum_{j=n+1-k}^n c_j$ for symbolical simplicity.
Then
$$
\mC_{\law\left(\sum_{j=1}^n c_j X_{t_j}\right)}(z)
= \mC_{\law\left(\sum_{k=1}^n c_{n,k}\left(X_{t_{n+1-k}}-X_{t_{n-k}}\right)\right)}(z)
, \q z \in \C^{-}.
$$
Since $\xt$ is freely additive, the above is
\begin{eqnarray*}
&=& \mC_{\bp_{k=1}^n D_{c_{n,k}} \law\left(X_{t_{n+1-k}}-X_{t_{n-k}}\right)}(z)\\
&=& \sum_{k=1}^n \mC_{\law\left(X_{t_{n+1-k}}-X_{t_{n-k}}\right)}(c_{n,k}z)\\
&=& \sum_{k=1}^n \mC_{\law\left(X'_{t_{n+1-k}}-X'_{t_{n-k}}\right)}(c_{n,k}z)\\
&=& \sum_{k=1}^n \mC_{D_{c_{n,k}}\law\left(X'_{t_{n+1-k}}-X'_{t_{n-k}}\right)}(z)\\
&=&  \mC_{\bp_{k=1}^n D_{c_{n,k}}\law\left(X'_{t_{n+1-k}}-X'_{t_{n-k}}\right)}(z)\\
&=&  \mC_{\law\left(\sum_{k=1}^n c_{n,k}(X'_{t_{n+1-k}}-X'_{t_{n-k}})\right)}(z)\\
&=& \mC _{\law\left( \sum_{j=1}^n c_jX'_{t_j}\right)}(z), \q\q z \in \C^{-}.
\end{eqnarray*}
\qed

The equivalence of the two notions of selfsimilarity for free additive processes can be given
as follows.

\begin{thm}
Suppose that $\xt$ is a free additive process and marginally $H$-selfsimilar such that $\law(X_{at}) = \law(a^HX_t)$ for all $a>0$ and $t\ge 0$.
Then $\xt$ is $H$-selfsimilar in the sense of Definition~\ref{free-ss}.
\end{thm}

{\it Proof.}
It is enough to take $\{X_{at}, t\ge 0\}$ and $\{a^HX_t, t\ge 0\}$ as $\xt$ and 
$\{X'_t, t\ge 0\}$ in Theorem~\ref{same}, respectively.
\qed

\vskip 3mm
We finish this section with some examples of selfsimilar non-commutative 
stochastic processes.

\begin{defn}
A $\mu\in\mathcal P (\R)$ is called a free stable distribution if for any $a>0$ and $b>0$, there
exists $c>0$ depending on $a$ and $b$ such that $D_a\mu \bp D_b\mu = D_c\mu$.
\end{defn}

\begin{rem}
There exists a unique $\al\in (0,2]$ such that $c$ above can be expressed as
$c= (a^{\al}+b^{\al})^{1/\al}$.
In this case, we call $\al$ the index of free stability and the free stable distribution $\mu$ to be free $\al$-stable.
When $\al=2$, if the mean is $0$ and the variance is $1$, then $\mu$ is a standard semicircle distribution such that $\mu(dx) = \frac1{2\pi}\sqrt{4-x^2}{\bf 1}_{\{|x|\le 2\}}dx$
 and the free $2$-stable L\'evy process is called a free Brownian motion.
\end{rem}

\begin{ex}
A free Brownian motion is $\frac12$-selfsimilar.
\end{ex}

\begin{ex}
A free $\al$-stable L\'evy process is $\frac1{\al}$-selfsimilar.
\end{ex}

\begin{ex}
(See Definition 2.2, Theorem 2.3 and Theorem 2.6 of \cite{G06}.)
Let $\xt$ be a non-commutative stochastic process.
It is called a standard semicircle process if for any $n\in\N, c_1,... , c_n\in\R,
0\le t_1< \cdots < t_n, \sum_{j=1}^n c_jX_{t_j}$ has a standard semicircle distribution.
A standard semicircle process is called a fractional free Brownian motion with index $H\in (0,1)$, if
$$
\tau (X_tX_s) = \frac12 \left( t^{2H} + s^{2H} - |t-s|^{2H}\right). \q t,s\ge 0.
$$
It is $H$-selfsimilar. When $H\not=\frac12$, it is not a free additive process, but when $H=\frac12$, it is a free additive process, which is nothing but a free Brownian motion.
\end{ex}

\begin{ex}
In the paper by Nourdin and Taqqu \cite{NT14}, they discuss not only the fractional free Brownian motion,
which is a selfsimilar semicircle process, but also the non-commutative Tchbycheff processes,
which are selfsimilar non-semicircle processes, as limiting processes of
so-called non-central limit theorems in a free probability setting.
Readers who want to know the details are suggested to read \cite{NT14}.
\end{ex}

In classical probability theory, the study of selfsimilar processes has a long history.
(See e.g. Chapter 7 of \cite{ST94} and \cite{EM02}.)
There are many examples of selfsimilar processes and problems to which the marginal
selfsimilarity is not enough, and it is a reason why we have introduced Definition~\ref{free-ss}.


\vskip 10mm
\section{A relation between selfsimilar free additive processes and freely selfdecomposable distributions}

We now consider freely selfdecomposable distributions. 

\begin{defn}\label{defSD}
Let $\mu \in \mathcal P (\R)$.
We say that $\mu$ is freely selfdecomposable (also called $\boxplus$-selfdecomposable), 
if for any $c\in (0,1)$ there exists a
unique $\rho_c\in\mathcal P(\R)$ such that
$$
\mu = D_c\mu\boxplus \rho_c.
$$
\end{defn}

\begin{rem}
(Proposition 4.22 (ii) of \cite{BT06})
Any free stable distribution mentioned in the previous section is freely selfdecomposable.
\end{rem}

Denote by $L(\bp)$ the class of all freely selfdecomposable distributions on $\R$.
Definition \ref{defSD} can be characterized in terms of the free cumulant transform  $\mC_{\mu}$ as follows.

\begin{prop}[Remark 4.3 of  \cite {BT02}]
$\mu \in \mathcal P (\R)$ is $\bp$-selfdecomposable if and only if for any 
$c\in (0,1)$, there exists a unique $\rho_c\in\mathcal P(\R)$ such that
$$
\mC_{\mu} (z) = \mC_{\mu}(cz) + \mC_{\rho_c}(z), \q z\in\C^- .
$$
\end{prop}

\begin{prop}[Theorem 3.29 of  \cite{F06}]\label{free-self}
If a non-commutative stochastic process $\{X_t, t\ge0\}$ of selfadjoint operators, affiliated with some 
$W^*$-probability space $(\mA, \tau)$, is marginally selfsimilar in the sense of \eqref{mss} by Fan \cite{F06} and freely additive,
then for each $t\ge 0$, $\law(X_t)$ is $\bp$-selfdecomposable.
\end{prop}


Theorem 16.1 (i) of Sato \cite{S99} is a counterpart of Proposition~\ref{free-self} 
above in classical probability theory.
In his Theorem 16.1 (ii), Sato \cite{S99} also shows the converse of (i).
Here we show the converse of Proposition \ref{free-self} for selfsimilar free additive processes.

\begin{thm}\label{main}
If $\mu\in\mathcal P(\R)$ is a non-trivial freely selfdecomposable distribution, then for any $H>0$,
there exists, uniquely in law, a non-trivial $H$-selfsimilar free additive process 
$\{X_t, t\ge 0\}$ of selfadjoint operators affiliated with some $W^*$-probability
space $(\mA,\tau)$ such that $\law(X_1)=\mu$.
\end{thm}

For the proof, we need the following.

\begin{prop}[\cite{B98},  Remark 4.28 of \cite{BT06}]\label{B98}\label{const}
Given a family of probability measures $\{\mu_t\,\, |\,\, t\ge 0\} \cup \{\mu_{s,t}
\,\, | \,\, 0\le s<t\}$ on $\R$ 
satisfying
$$
\mu_0=\dl_0, \q \mu_{0,t}=\mu_t, \q \mu_{t,t}=\mu_0,
$$
\begin{align}\label{1st}
\mu_s\bp\mu_{s,t} = \mu_t, \q \mu_{s,r}\bp\mu_{r,t}=\mu_{s,t}\q  {\rm for}\quad  0<s<r<t,
\end{align}
there exists a free additive process $\{X_t, t\ge 0\}$ affiliated with a
$W^*$-probability space $(\mathcal{A}, \tau)$ such that $\mu_t=\law(X_t)$
and $\mu_{s,t}=\law(X_t-X_s)$ for all  $0\le s\le t$.
\end{prop}

Biane \cite{B98} mentions that this follows from the free product construction and an inductive limit argument. 

\vskip 3mm

{\it Proof of Theorem~\ref{main}.}
(Almost the same argument as for the case of classical probability theory given in 
Theorem 16.1 (ii) of Sato \cite{S99}.)

Step 1. (Construction of a free additive process by Proposition \ref{const}.)

Let $\mu$ be a non-trivial $\bp$-selfdecomposable distribution and $H>0$.
Then $\mu$ is $\bp$-infinitely divisible,  (see Proposition 4.26 of \cite{BT06}).
Thus for any $c\in (0,1)$, there exists a unique probability measure $\rho _c$ such that
\begin{eqnarray}\label{Cselfdec}
\mC_{\mu}(z) = \mC_{\mu}(cz) + \mC_{\rho_c}(z), \q\q z\in \C^-.
\end{eqnarray}
It follows that $\rho_c$ is continuous in $c\in (0,1)$ in the sense of weak convergence.
It comes from a similar discussion with Lemma 4.24 in \cite{BT06}.

In the following, we use two properties of $\mC _{\mu}(z)$, which are already
mentioned in Section 2, frequently.
$$
\mC _{\mu_1\bp\mu_2}(z) = \mC_{\mu_1}(z) + \mC_{\mu_2}(z), \q z \in \C^{-},
$$
and
$$
\mC_{D_c\mu}(z) = \mC_{\mu}(cz), \q z \in \C^{-}.
$$
Define, for $t>0$ and $0\le s<t$, $\mu_t$ and $\mu_{s,t}$ by
$$
\mC_{\mu_t}(z) = \mC_{\mu}(t^Hz), \q z \in \C^{-},
$$
and
\begin{eqnarray}\label{43}
\mC_{\mu_{s,t}}(z) = \mC_{\rho_{{}_{(s/t)^H}}}(t^Hz), \q z \in \C^{-},
\end{eqnarray}
where $\rho_{(s/t)^H}$ is the one in \eqref{Cselfdec}.
Then
\begin{eqnarray}
\mC_{\mu_t}(z) & = & \mC _{\mu}\left( \left(\frac{s}{t}\right)^Ht^Hz\right) + 
\mC_{\rho_{{}_{(s/t)^H}}}(t^Hz)\nonumber\\
& = & \mC_{\mu}(s^Hz) + \mC_{\mu_{s,t}}(z)\nonumber\\
&=& \mC_{\mu_s}(z) + \mC_{\mu_{s,t}}(z), \q z \in \C^{-}.  \label{2nd}
\end{eqnarray}
Further define $\mu_0=\dl_0, \mu_{0,t}=\mu_t$ and $\mu_{t,t}=\dl_0$.
Then it follows from \eqref{2nd} that
\begin{align}\label{3rd}
\mu_t = \mu_s\bp\mu_{s,t}
\end{align}
for $0\le s\le t$ and $\mu_t$ is continuous in $t\ge 0$.

By Lemma 4.24 of \cite{BT06}, we see that if $c_n\to 1$ as $n\to \infty$, then $\mu_{c_n}
\to \dl_0$ as $n\to \infty$.
Thus,
$$
\mC_{\mu_{s,t}}(z) \to \mC_{\dl_0}(z), \q z \in \C^{-}
\q {\rm as}\q s\uparrow t \q  \text{or} \q t\downarrow s
$$
by \eqref{Cselfdec} and \eqref{43}, and we have
$$
\mu_{s,t}\to \dl_0, \q\text{as}\q s\uparrow t\q  \text{or} \q t\downarrow s.
$$
For  $s\le r\le t$, we have, from \eqref{2nd},
\begin{eqnarray}
\mC_{\mu_{s,t}}(z) &=& \mC_{\mu_t}(z) - \mC_{\mu_s}(z) \nonumber \\
&=& \mC_{\mu_t}(z) - \mC_{\mu_r}(z) + \mC_{\mu_r}(z) - \mC_{\mu_s}(z)\nonumber \\
&=& \mC_{\mu_{s,r}}(z) + \mC_{\mu_{r,t}}(z), \q\q z \in \C^{-}, \nonumber
\end{eqnarray}
meaning
\begin{align}\label{4th}
\mu_{s,t} = \mu_{s,r}\bp \mu_{r,t}.
\end{align}
Equations \eqref{3rd} and \eqref{4th} are \eqref{1st}, and thus by Proposition~\ref{B98}, there exists, uniquely in law,
a free additive process $\{X_t, t\ge 0\}$ such that $\law(X_1)= \mu_1$ and $\law(X_t-X_s) = \mu_{s,t}$.

Step 2. (Proof of the selfsimilarity.)

By the definition of $\mu_t$, we have, for $a>0$,
$$
\mC_{\mu_{at}}(z) = \mC_{\mu}(a^Ht^Hz) = \mC_{D_{a^H\mu_t}}(z), \q z \in \C^{-},
$$
namely, $\law(X_{at}) = \law(a^HX_t)$ for all $t\ge 0$, meaning the marginal
selfsimilarity of $\xt$.

Hence $\xt$ is $H$-selfsimilar.
\qed


\vskip 10mm
\section{Stochastic integrals with respect to free additive processes}

In classical probability theory, Sato \cite{S04} defined stochastic integrals with respect to
(classical) additive processes $\{Y_t\}$ rigorously.
He pointed out that some extra mild condition is needed to define the integrals.
The condition was named ``natural'' by him, meaning that the location parameter 
$\gm_t$ in the characteristic triplet of the infinitely divisible distribution $\law(Y_t)$
is locally of bounded variation in $t$.
He also showed that $\{Y_t\}$ is natural if and only if $\{Y_t\}$ is 
martingale,
and that selfsimilar additive processes are natural.

We here use the same concept for the free additive process, namely,
a free additive process $\{X_t\}$ is called natural if the location parameter $\gm_t$ in the
free characteristic triplet of $\law(X_t)$ is locally of bounded variation in $t$.

We now introduce the Bercovici-Pata bijection, which will be used in this section
for the proof of many statements.

\begin{defn}[\cite{BP90}]
By the Bercovici-Pata bijection $\Ld : ID(*) \to ID(\bp)$, we define the mapping as follows:
Let $\mu$ be a probability measure in $ID(*)$ with the generating pair $(\gm, \sg)$
and the characteristic triplet $(a,\nu,\eta)$.
Then $\Ld(\mu)$ is the probability measure in $ID(\bp)$ with the free generating pair $(\gm, \sg)$ and the free characteristic triplet $(a,\nu,\eta)$.
\end{defn}

By using the Bercovici-Pata bijection, it can be seen that a selfsimilar free additive process is also natural.

We collect here several properties of $\Ld$ needed later.
All below are from \cite{BT06}.

\begin{thm}\label{BerPata}
The Bercovici-Pata bijection $\Ld$ has the following properties.
\begin{enumerate}[\rm (i)]
\item 
If $\mu_1,\mu_2\in ID(*)$, then $\Ld(\mu_1*\mu_2) = \Ld(\mu_1)\bp\Ld(\mu_2)$.
\item 
If $\mu\in ID(*)$ and $c\in\R$, then $\Ld(D_c\mu)=D_c\Ld(\mu)$.
\item 
For any constant $c\in\R$, we have $\Ld(\dl_c)=\dl_c$.
\item
The bijection $\Ld$ is invariant under affine transformations, i.e. if $\mu\in ID(*)$ and 
$\psi:\R\to\R$ is an affine transformation, then
$$
\Ld(\psi(\mu))= \psi(\Ld(\mu)).
$$
\item
The bijection $\Ld$ is a homeomorphism with respect to weak convergence.
In other words, if $\mu$ is a measure in $ID(*)$ and $\{\mu_n\}$ is a sequence of measures in 
$ID(*)$, then $\mu_n\to\mu$ weakly as $n\to\infty$ if and only if 
$\Ld(\mu_n) \to \Ld(\mu)$ weakly as $n\to \infty$.
\end{enumerate}
\end{thm}

As for the correspondence between classical and free additive processes through the Bercovici-Pata bijection,
the following is also true.
Here, as mentioned in Remark~\ref{additive-free}, we use the $\bp$-infinite divisibility 
of marginal distributions of free additive processes.
The proof can be done by the same argument as for Proposition 5.15 of \cite{BT06}
for free L\'evy processes.

\begin{prop}\label{bijec}
Let $\xt$ be a free additive process affiliated with a $W^*$-probability space $(\mA,\tau)$.
Then there exists a (classical) additive process $\yt$ with $\law (Y_t) = \Ld^{-1}(\law(X_t))$ for each $t\ge 0$.
Conversely, for a (classical) additive process $\yt$, there exists a free additive process
$\xt$ with $\law(X_t) = \Ld (\law(Y_t))$ for each $t\ge 0$.
\end{prop}

In addition to Proposition~\ref{bijec}, we have also the following.
Here we use the $\bp$-infinite divisibility of increments of free additive processes,
shown in Theorem~\ref{dif-id}.

\begin{prop}\label{bijec2}
Proposition~\ref{bijec} also holds with the replacement of $\law (Y_t) = \Ld^{-1}(\law(X_t))$
and  $\law(X_t) = \Ld (\law(Y_t))$ by $\law (Y_t-Y_s) = \Ld^{-1}(\law(X_t-X_s))$ and
$\law(X_t-X_s) = \Ld (\law(Y_t-Y_s))$, respectively, where $0\le s<t$.
\end{prop}

{\it Proof.}
Note that, in terms of cumulant functions, $C_{\mu}(\theta) =\log \widehat{\mu}(\theta), \theta\in\R$,
in classical probability theory,
$$
C_{\Ld^{-1}(\law(X_s))} (\theta)+ C_{\Ld^{-1}(\law(X_t-X_s))}(\theta) = C_{\Ld^{-1}(\law(X_t))}(\theta), \q \theta\in\R,
$$
and
$$
C_{\law(Y_s)} (\theta)+ C_{\law(Y_t-Y_s)}(\theta) = C_{\law(Y_t)}(\theta), \q 
\theta\in\R.
$$
Since $\law(Y_t) = \Ld^{-1}(\law(X_t))$
by Proposition~\ref{bijec},
we have $\law(Y_t-Y_s) = \Ld^{-1}(\law(X_t-X_s)).$
\qed


\vskip 3mm
Let $\xt$ be a natural free additive process and let $f : [A,B] \to \R$ be a continuous function defined on an interval $[A,B]\subset [0,\infty)$.
We are going to define the stochastic integral $\int_A^Bf(t)dX_t$ as the limit of approximating Riemann sums with the help of the Bercovici-Pata bijection.

\begin{thm}\label{stint}
Let $\xt$ be a natural free additive process and $f$ an $\R$-valued measurable function
defined on an interval $[A,B]$, $0\le A < B <\infty$. 
Let $n\in\N$ and let $\mathcal D_n=\{t_{n,0}, t_{n,1}, ..., t_{n,n}\}$ be a subdivision of the
interval $[A,B]$, namely,
$$
A= t_{n,0}<t_{n,1}<\cdots < t_{n,n}=B.
$$
Assume that
$$
\lim_{n\to\infty} \max_{1\le j\le n}(t_{n, j}-t_{n, j-1})=0.
$$
Furthermore, for each $n$, choose intermediate points
$$
t_{n,j}^{\#}\in (t_{n,j-1}, t_{n,j}], \q j=1,2,...,n.
$$
Then there exists uniquely a non-commutative random variable $I$ with its law in $ID(\bp)$
to which the corresponding Riemann sums
$$
I_n:= \sum_{j=1}^n f(t_{n,j}^{\#})(X_{t_{n,j}}-X_{t_{n,j-1}})
$$
converge in probability as $n\to\infty$.
\end{thm}

We call $I$ above the stochastic integral of $f$ with respect to $\xt$ and denote it
by $\int_A^Bf(t)dX_t$.
In the proof below, we also use the notation
$$
*_{j=1}^n\mu_j = \mu_1* \cdots *\mu_n
$$
for $\mu_1,...,\mu_n\in\mathcal P (\R)$.

\vskip 3mm
{\it Proof of Theorem~\ref{stint}.}
For any $n,m\in\N$, we consider the subdivision
$$
A=s_0<s_1<\cdots <s_{p(n,m)}=B,
$$
which consists of the points in $\mathcal D_n\cup\mathcal  D_m$,
and for each $j\in\{1,2,...., p(n,m)\}$, we choose $s_{n,j}^{\#}\in \{t_{n,k}^{\#} | k=1,2,..., n\}$
and $s_{m,j}^{\#}\in \{t_{m,k}^{\#} | k=1,2,..., m\}$
such that 
$$
\sum_{j=1}^{p(n,m)} f(s_{n,j}^{\#})(X_{s_j}-X_{s_{j-1}})
= \sum_{k=1}^m f(t_{m,k}^{\#})(X_{t_{n,k}} - X_{t_{n,k-1}}).
$$
Put
$$
J_n:= \sum_{j=1}^{p(n,m)} f(s_{n,j}^{\#}) (X_{s_j}- X_{s_{j-1}})
\q{\rm and}\q J_m:= \sum_{j=1}^{p(n,m)} f(s_{m,j}^{\#}) (X_{s_j}- X_{s_{j-1}}).
$$
Then
$$
J_n-J_m = \sum_{j=1}^{p(n,m)} \left (f(s_{n,j}^{\#}) - f(s_{m,j}^{\#})\right )(X_{s_j}-X_{s_{j-1}}).
$$
Consider a natural (classical) additive process $\yt$ such that $\law( Y_t) = \Ld^{-1}(\law(X_t))$,
which is possible by Proposition~\ref{bijec}.
For each $n$, we form the Riemann sum
$$
K_n = \sum _{j=1}^n f(s_{n,j}^{\#})(Y_{t_{n,j}}-Y_{t_{n,j-1}}).
$$
Then, for any $n,m\in\N$, we have also
$$
K_n - K_m = \sum_{j=1}^{p(n,m)} \left( f(s_{n,j}^{\#}) -f(s_{m,j}^{\#})\right) (Y_{s_j}- Y_{s_{j-1}}).
$$
Thus
$$
\law(K_n-K_m) = *_{j=1}^{p(n,m)} D _{f(s_{n,j}^{\#}) -f(s_{m,j}^{\#})}\law(Y_{s_j}-Y_{s_{j-1}})
$$
and therefore
\begin{eqnarray*}
\Ld\left( \law(K_n-K_m)\right)&=&
\bp _{j=1}^{p(n,m)} D _{f(s_{n,j}^{\#}) -f(s_{m,j}^{\#})}\Ld \left( \law(Y_{s_j}-Y_{s_{j-1}})\right )\\
&=& \law\left( \sum_{j=1}^{p(n,m)}(f(s_{n,m}^{\#}) - f(s_{m,j}^{\#})) (X_{s_j}-X_{s_{j-1}})\right)\\
&=& \law(J_n-J_m),
\end{eqnarray*}
where we have used Proposition~\ref{bijec2}.
It can be shown by the facts in \cite{S04} that
$$
\law(K_n-K_m) \to \delta _0 \q {\rm as}\,\, n,m\to \infty.
$$
Then by the same argument as in the proof of Theorem 6.1 in \cite{BT06},
there exists a selfadjoint operator $I$ in $\bar{\mA}$ such that $I_n \to I$ in probability.

The fact that the operator $I$ above does not depend on the choice of subdivision 
$\mathcal D_n$ or intermediate points $\{ t_{n,j}^{\#}\}$ can be shown again by the same reasoning as in the proof of Proposition 6.1 in \cite{BT06}.
The free infinite divisibility of $I = \int_A^B f(t)dX_t$ is a consequence of the closedness of $ID(\bp)$ with respect to weak convergence.
The proof is thus completed.
\qed

\begin{cor}\label{bijec3}
Let $\yt$ be a natural (classical) additive process and let $\xt$ be a corresponding
natural free additive process with marginal distributions $\Ld (\law(Y_t))$.
Then for any compact interval $[A,B]\subset [0,\infty)$ and any continuous function
$f : [A,B]\to \R$,
$$
\law\left( \int_A^B f(t)dX_t\right ) = \Ld \left ( \law\left( \int_A^B f(t)dY_t\right)\right).
$$
\end{cor}

The proof can be carried out in the same manner as in the proof of Corollary 6.2
of \cite{BT06}.


\vskip 10mm
\section{Background driving free L\'evy processes of freely selfdecomposable distributions}

\subsection{Construction of background driving free L\'evy processes of freely selfdecomposable distributions}

Let $\mu\in ID(\bp)$
and let $(\gm,\sg)$ be the free generating pair of $\mu$ and $(a,\nu,\eta)$ the free characteristic triplet of $\mu$, respectively.
Then
$$
\int_{\R\setminus [-1,1]}\log (1+|x|)\sg(dx) <\infty
$$
if and only if 
$$
\int_{\R\setminus [-1,1]}\log (1+|x|)\nu(dx) <\infty,
$$
because 
$$
\sg(dx) = s\dl_0(dx) + \frac{x^2}{1+x^2}\nu(dx).
$$

We start with the following statement.

\begin{prop}[Theorem 6.5 of \cite{BT06}]
Let $X$ be a selfadjoint operator affiliated with a $W^*$-probability space $(\mA, \tau)$.
Then $\law(X)$ is $\bp$-selfdecomposable if and only if $X$ has a representation in law
in the form 
$$
\law(X) = \law\left( \int_0^{\infty} e^{-t}dZ_t\right)
$$
for some free L\'evy process $\zt$ affiliated with some $W^*$-probability space
satisfying 
\begin{align}\label{logmoment}
\int_{\R\setminus [-1,1]}\log (1+|x|)\nu (dx) < \infty,
\end{align}
where 
$\nu$ is the free L\'evy measure of $\law(Z_1)$.
\end{prop}

We call the free L\'evy process $\zt$ the background driving free L\'evy process 
 of $\law(X)\in L(\bp)$ as in classical probability theory.

In classical probability theory, it was shown that the background driving L\'evy process of
(a classical random variable) $X$ can be constructed as follows.  

\begin{prop}
{\rm (a) (}\cite{S91}, also Theorem 16.1 of \cite{S99}.{\rm )}
$\law(X)$ is non-trivial selfdecomposable if and only if for any $H>0$, there exists, 
uniquely in law, a non-trivial $H$-selfsimilar additive process $\xt$ such that
$\law(X_1)=\law(X)$.

\n
{\rm (b) (}See \cite{JPY02} and \cite{MS03}.{\rm)}
Let
\begin{equation}\label{Zt}
Z_t = \int_1^{e^t} u^{-H}dX_u,\q t\in\R.
\end{equation}
Then

\n
{\rm (b1)} \q $\{ Z_t, t\in\R\}$ is a L\'evy process.

\n
{\rm (b2)} \q $\int_{\R\setminus [-1,1]}\log (1+|x|)\nu (dx) < \infty$, where the measure
$\nu$ is the L\'evy measure of the L\'evy process $\{Z_t\}$ shown in {\rm (b1)}.

\n
{\rm (b3)} \q $ X_t = \int_{-\infty}^{\log t} e^{uH}dZ_u, \,\, t> 0,$
and $X=\int_0^{\infty}e^{-uH}dZ_u$.

\n
Thus, $\zt$ in \eqref{Zt} is the background driving L\'evy process of $\law(X)\in L(*)$.
\end{prop}

We now study the same problem in free probability theory.
In the following, the symbol $\stackrel{\rm law}{=}$ means the equality in law.
The free version of the part (a) above has already been shown in Proposition~\ref{free-self} and Theorem~\ref{main}.
So, we prove the free version of the part (b) above as follows.

\begin{thm}\label{main2}
Let $X$ be a selfadjoint operator affiliated with a $W^*$-probability space $(\mA, \tau)$.
Suppose $\law(X)$ is $\bp$-selfdecomposable and $\xt$ is the corresponding 
$H$-selfsimilar free additive process such that $\law(X_1)=\law(X)$.
Let
\begin{equation}\label{Ztfree}
Z_t = \int_1^{e^t} u^{-H}dX_u,\q t\in\R.
\end{equation}
Then

\n
{\rm (b1)} \q $\{ Z_t, t\in\R\}$ is a free L\'evy process.

\n
{\rm (b2)} \q $\int_{\R\setminus [-1,1]}\log (1+|x|)\nu (dx) < \infty$, where the measure
$\nu$ is the free L\'evy measure of the free L\'evy process $\{Z_t\}$ shown in {\rm (b1)}.

\n
{\rm (b3)}  
\begin{eqnarray}\label{b3}
X_t = \plimR\int_{-R}^{\log t}e^{uH}dZ_u, \,\, t>0,
\end{eqnarray} 
and 
$$
X\stackrel{\rm law}{=}\plimR \int_0^R u^{-uH}dZ_u =: \int_0^{\infty}e^{-uH}dZ_u,
$$
where $\plim$ means the limit in probability.
Thus, $\{ Z_t\}$ in \eqref{Ztfree} is the background driving free L\'evy process of $\law(X)\in L(\bp)$.
\end{thm}

We prepare three lemmas.

\begin{lem}\label{64}
The following change of variables is possible: for $-\infty <s<t<\infty$,
$$
\int_{e^s}^{e^t}u^{-H}dX_u 
\stackrel{\rm law}{=}
 \int_1^{e^{t-s}} (e^sv)^{-H}d_v\left(X_{e^sv}\right).
$$
\end{lem}

{\it Proof.}
The statement is easily verified by noticing that the integrals are defined as limits of 
Riemann sums.
Actually, for $n\in \N$ and $e^s=u_{n,0}< u_{n,1}<\cdots <u_{n,n} =e^t$ satisfying
$\lim_{n\to\infty} \max_{1\le j \le n} (u_{n,j} -u_{n,j-1})=0$,
\begin{eqnarray*}
\int_{e^s}^{e^t}u^{-H}dX_u &\stackrel{\rm law}{=}& \plim_{n\to\infty}\sum_{j=1}^nu_{n,j}^{-H} (X_{u_{n,j}}  - X_{u_{n,j-1}})\\
&=& \plim_{n\to\infty}\sum_{j=1}^n (e^sv_{n,j})^{-H}\left( X_{e^sv_{n,j}}-X_{e^sv_{n,j-1}}\right)\\
&\stackrel{\rm law}{=}& \int_1^{e^{t-s}}(e^sv)^{-H}d_v\left(X_{e^sv}\right).
\end{eqnarray*}
\qed

\begin{lem}\label{H-ss}
We have
$$
\int_1^{e^{t-s}}(e^sv)^{-H}d_v\left(X_{e^sv}\right) \stackrel{\rm law}{=}
\int_1^{e^{t-s}}v^{-H}dX_v = Z_{t-s}.
$$
\end{lem}

{\it Proof.}
As the same as before, the statement is easily verified by considering the integrals as 
limits of Riemann sums.
Actually, for $n\in \N$ and $1=v_{n,0}< v_{n,1}<\cdots <v_{n,n} =e^{t-s}$ satisfying
$\lim_{n\to\infty} \max_{1\le j \le n} (v_{n,j} -v_{n,j-1})=0$,
\begin{eqnarray}\label{free-ss1}
\int_1^{e^{t-s}}(e^sv)^{-H}d_v\left(X_{e^sv}\right) \stackrel{\rm law}{=} \plim_{n\to\infty}\sum_{j=1}^n
(e^sv_{n,j})^{-H}\left(X_{e^sv_{n,j}}-X_{e^sv_{n,j-1}} \right).
\end{eqnarray}
By $H$-selfsimilarity of $\{X_t\}$, the above is equal in law to
\begin{eqnarray}\label{free-ss2}
 \plim_{n\to\infty}\sum_{j=1}^n v_{n,j}^{-H}\left(X_{v_{n,j}} - X_{v_{n,j-1}}\right)
 {=} \int_1^{e^{t-s}}v^HdX_v = Z_{t-s}.
\end{eqnarray}
\qed

\vskip 3mm
{\it Proof of Theorem~\ref{main2} (b1).}
The increments of $\{Z_t, t\in\n\R\}$ are integrals with respect to a free additive process
$\{X_t, t\in\R\}$ on disjoint intervals.
Thus $\{Z_t, t\in\R\}$ has the freely independent increment property.
For $-\infty < s<t<\infty$, we have
$$
Z_t-Z_s = \int_{e^s}^{e^t} u^{-H}dX_u 
\stackrel{\rm law}{=}
 \int_1^{e^{(t-s)}} (e^sv)^{-H}d_v(X_{e^sv}),
$$
by Lemma~\ref{64}.
Furthermore, by Lemma~\ref{H-ss},
the above is equal in law to
$$
\int_1^{e^{t-s}}v^{-H}dX_v = Z_{t-s}.
$$
Hence,  $\{Z_t, t\in\R\}$ has the stationary increment property. 
It is trivial that $Z_0=0$, and we also see that $\{Z_t, t\in\R\}$ is stochastically continuous, which can be shown by the continuity of the Bercovici-Pata bijection as used in the proof of Theorem 6.26 of \cite{BT06}. 
Thus,  $\{Z_t, t\in\R\}$ is a (two-sided) L\'evy process.
\qed

\begin{lem}\label{relation}
Let $\{Z_t,t\in\R\}$ be the one defined in \eqref{Ztfree}.
Then, for $0<t_0<t$,
$$
\int_{\log t_0}^{\log t} e^{uH}dZ_u  \stackrel{\rm law}{=}
X_t -X_{t_0}.
$$
\end{lem}

{\it Proof.}
Let $\{Y_t\}$ be a natural (classical) additive process with marginal distributions 
$\law(Y_t)= \Ld^{-1}(\law(X_t))$ and let
$$
W_t= \int_1^{e^t}u^{-H}dY_u, \q t\in\R.
$$
Then it is known from the proof of (1.6) of Theorem 1.1 in \cite{MS03} that, for $0<t_0<t$,
$$
\int_{\log t_0}^{\log t}e^{uH}dW_u = Y_t - Y_{t_0}, \,\, {\rm almost \,\,surely.}
$$
Then by the Bercovici-Pata bijection, we have $\Ld(\law(W_t))=\law(Z_t)$ and $\Ld(\law(Y_t))=\law(X_t)$,
and thus
$$
\Ld\left(\law\left(\int_{\log t_0}^{\log t}e^{uH}dW_u\right)\right) 
= \Ld(\law(Y_t - Y_{t_0})),
$$
implying the conclusion of this lemma
$$
\law\left(\int_{\log t_0}^{\log t} e^{uH}dZ_u\right) =
\law(X_t -X_{t_0}),
$$
where we have used Proposition~\ref{bijec2} and Corollary~\ref{bijec3}.
\qed

\vskip 3mm
{\it Proof of Theorem~\ref{main2} (b3) and (b2).}

(b3).
By Lemma~\ref{relation}, for $0<t_0<t$, we have
$$
\int_{\log t_0}^{\log t} e^{uH}dZ_u
\stackrel{\rm law}{=}
 X_t -X_{t_0}.
$$
Since, $\lim_{t_0{\downarrow}0}X_{t_0}=0$ in probability,
$$
\int_{-R}^{\log t} e^{uH}dZ_u
$$
converges in probability as $R\to\infty$.
Thus,
$$
X_t = \plimR\int_{-R}^{\log t}e^{uH}dZ_u ,
$$
which is \eqref{b3}.
Thus, by taking $t=1$, we have
$$
X \stackrel{\rm law}{=} X_1 = \int_{-\infty}^0 e^{uH}dZ_u
\stackrel{\rm law}{=} \int_0^{\infty} e^{-vH}dZ_v.
$$

(b2).
By Proposition 6.4 of \cite{BT06}, it follows from (b3) that
$$
\int_{\R\setminus [-1,1]}\log (1+|x|)\nu (dx) < \infty.
$$
\qed

\vskip 3mm
The proof of Theorem~\ref{main2} is now completed. 


\subsection{The free cumulant transforms of the background driving free L\'evy processes of freely selfdecomposable distributions}

In classical probability theory, Jurek \cite{Jurek01} gave a new description of the characteristic function of the background driving L\'evy processes of the selfdecomposable distributions as follows:   Let $\mu \in\mathcal{P}(\R)$ be a selfdecomposable distribution and let $\rho \in\mathcal{P}(\R)$ be the distribution of its background driving L\'evy process. Then 
\begin{align*}
\wh{\rho}(\theta) = \exp\left\{\theta \frac{\wh{\mu}'(\theta)}{\wh{\mu}(\theta)}\right\},\quad \theta \in\R,
\end{align*}
namely,
$$
\log \widehat{\rho}(\theta) = \theta \frac{d}{d\theta}\log \widehat{\mu}(\theta), \q \theta\in\R.
$$
Here we show a similar relation in the free probability setting, which proves to be helpful in discussing examples in the next section.
Note that the free cumulant of the background driving free L\'evy process $\zt$ satisfies
$\mC_{\law(Z_t)}(z) = t\mC_{\law(Z_1)}(z)$,
where the equality holds by the independent and stationary increment property of
the free L\'evy process $\zt$.
The following holds for $\mC_{\law(Z_1)}(z)$.

\begin{thm}\label{thmformula}
Under the notation in Proposition 6.1, we have 
\begin{align}\label{formula2}
\mC_{\law(Z_{1})}(z) 
= z\frac{d}{dz}\mC_{\law(X)}(z),
\quad z\in\C^{-}.
\end{align}
\end{thm}
\begin{rem}
It follows from Theorem~\ref{main2} that
$Z_{1}$ is freely infinitely divisible.
On the other hand, 
as mentioned after Definition~\ref{ana-ext}, by Theorem 5.10 of \cite{BV93}, $\mu\in ID(\boxplus)$ if and only if $\phi_{\mu}(z)$ has an analytic extension to $\C^{+}$ with values in $\C^{-}\cup\R$.
So, $\phi_{\law(Z_1)}(z) = z\mC_{\law(Z_{1})}(z^{-1}) = \frac{d}{dz}\mC_{\law(X)}(z^{-1})$ has to be in $\C^{-}\cup\R$.
It is the same as Theorem 2.7 (ii) in \cite{HST19}, which is a necessary and sufficient condition for $\law(X) \in L(\boxplus)$. 
\end{rem}

{\it Proof of Theorem \ref{thmformula}.}
In the setting of Proposition 6.1, we have
\begin{align*}
\mC_{\law(X)}(z) 
&=\mC_{\law(\int_{0}^{\infty} e^{-t} dZ_{t})}(z)\\
&=\mC_{\lim_{A\to\infty}\law(\int_{0}^{A} e^{-t} dZ_{t})}(z)\\ 
&=\lim_{A\to\infty}\mC_{\law(\int_{0}^{A} e^{-t} dZ_{t})}(z)\\  
&=\lim_{A\to\infty}\mC_{\lim_{n\to\infty}
\law(\sum_{j=1}^{n}e^{-j\frac{A}{n}} (Z_{j\frac{A}{n}} - Z_{(j-1)\frac{A}{n}})}(z)\\  
&=\lim_{A\to\infty}\lim_{n\to\infty}\sum_{j=1}^{n}
\mC_{\law(Z_{\frac{A}{n}})}(e^{-j\frac{A}{n}}z)\\  
&=\lim_{A\to\infty}\lim_{n\to\infty}\sum_{j=1}^{n}
\frac{A}{n} \mC_{\law(Z_{1})}(e^{-j\frac{A}{n}}z)\\
&=\lim_{A\to\infty}\int_{0}^{A}
\mC_{\law(Z_{1})}(e^{-t}z)dt, 
\quad z\in\C^{-}.
\end{align*}
Let $(a,\nu,\eta)$ be the free characteristic triplet of $\law(Z_{1})$.
Recall that our L\'evy measure $\nu$ satisfies \eqref{logmoment}.
To show the existence of the limit above,
it is enough  to show the finiteness of
\begin{align*}
\int_{0}^{\infty}\left|\mathcal{C}_{\mathcal{L}(Z_{1})}(e^{-t}z)\right|  dt
&\leq \int_{0}^{\infty} \left(\int_{\R} \left|\frac{1}{1-e^{-t}zx}-1 -e^{-t}zx{\bf 1}_{[-1,1]} (x) \right|\nu(dx) \right)dt\\
&\leq \int_{0}^{\infty}\Big( \int_{\R}\Big( \left|\frac{1}{1-e^{-t}zx}-1 -e^{-t}zx{\bf 1}_{[-e^{t},e^{t}]} (x)\right| \\
&\q\q +\left| e^{-t}zx{\bf 1}_{[-e^{t},e^{t}]} (x) -e^{-t}zx{\bf 1}_{[-1,1]} (x) \right|\Big)\nu(dx)\Big) dt\\
& =: I_1(z) +|z| I_2(z),\quad z\in\C^{-}.
\end{align*}
say.

In order to show that $I_1(z)<\infty$, we first show
$$
f(z):= \left|\frac{1}{1-e^{-t}zx}-1 -e^{-t}zx{\bf 1}_{[-e^{t},e^{t}]} (x)\right| 
\le C \left( |e^{-t}x|^2\land 1\right), \q {\rm for\,\, any}\,\, z\in\C^{-},
$$
where $C$ is a positive constant that can be chosen uniformly on any compact subset in $\C^{-}$.
Fix $z\in\C^{-}$.
Since $\C^{-}$ is an open set, there exists a compact set $K$ such that $z\in K$ and $K\subseteq \C^{-}$.
Fix such a $K$.
When $x\in [-e^t,e^t]$,
\begin{align*}
f(z)& = \left| \frac1{1-e^{-t}xz} - 1 - e^{-t}xz\right|
\le \sup_{z\in K} \frac{|e^{-t}x|^2|z|^2}{|1-e^{-t}xz|}\\
& \le \frac{M}{m} |e^{-t}x|^2 
 =: C_1 |e^{-t}x|^2,
\end{align*}
say,
where $M= \max_{z\in K} |z|^2 (<\infty)$, $m= \min_{x\in [-e^t, e^t], z\in K} |1- e^{-t}xz| (>0)$
and $C_1\in (0,\infty)$.
When $x\in\R \setminus [-e^t, e^t]$,
\begin{align*}
f(z) & = \left |\frac1{1-e^{-t}xz} - 1 \right |
 \le \sup _{z\in K} \frac{|e^{-t}xz|}{|1- e^{-t}xz|}
 = \sup_{z\in K} \frac{|z|}{|e^tx^{-1} -z|}\\
& \le \frac{\left(u^2_{\max} + v^2_{\max}\right)^{1/2}}{v_{\min}}
 =: C_2,
\end{align*}
say,
where $u_{\max} = \max \{ |\Re(z)| \,|\, z\in K\} (< \infty),
v_{\max} = \max \{ - \Im(z) | z\in K\} (< \infty),
v_{\min} = \min \{ - \Im (z) |z\in K\} (>0)$
and $C_2\in(0, \infty)$.
Thus,
\begin{align*}
I_1(z) & \le (C_1+C_2) \int_0^{\infty}\left( \int_{\R} \left( |e^{-t}x|^2\land 1\right)\nu(dx) \right)dt\\
& = (C_1+C_2) \int_0^{\infty} \left(\int_{\{|x|\le e^t\}}|e^{-t}x|^2\nu(dx)
+ \int_{\{|x|>e^t\}} \nu(dx)\right)dt\\
& =: (C_1+C_2) (I_{11}+I_{12}),
\end{align*}
say.
We have
\begin{align*}
I_{11} & = \int_0^{\infty} e^{-2t}
\left( \int_{\{|x|\le 1\}} |x|^2\nu(dx) + \int_{\{1< |x|\le e^t\}}|x|^2\nu(dx)\right) dt\\
& =: I_{111} + I_{112},
\end{align*}
say.
Here $I_{111}<\infty$, since $\int_{\{|x|\le 1\}}|x|^2\nu(dx) <\infty$ by
a property of the free L\'evy measure.
We also have, by the change of integrations,
\begin{align*}
I_{112} & = \int_0^{\infty} e^{-2t} \left(\int_{\{1<|x|\le e^t \}}|x|^2\nu(dx)\right)dt\\
& = \int_{\{|x|>1\}} |x|^2\left(\int_{\log |x|}^{\infty} e^{-2t}dt\right) \nu(dx)\\
& = \frac{1}{2}\int_{\{|x|>1\}}\nu(dx)
< \infty.
\end{align*}
As to $I_{12}$, 

\begin{align*}
I_{12}
&=\int_{0}^{\infty} \left(\int_{\{|x|>e^{t}\}}\nu(dx)\right)dt 
= \int_{\{|x|>1\}}\left(\int_{0}^{\log |x|} dt\right) \nu(dx)\\
&\le \int_{\{|x|>1\}}\log(1+ |x|) \nu(dx)<\infty,
\end{align*}
by \eqref{logmoment}.

Thus we conclude that $I_1(z)<\infty$.
As to $I_2$, we have
\begin{align*}
I_2 & =  \int_{\R} \left(\int_{0}^{\infty} e^{-t} |x| {\bf 1}_{[-e^{t},e^{t}]\backslash [-1,1]} (x) dt\right)\nu(dx)\\
& =\int_0^{\infty} e^{-t}\left( \int_{\R}  |x| {\bf 1}_{[-e^{t},e^{t}]\backslash [-1,1]} (x)\nu(dx)\right)dt\\
&= \int_0^{\infty} e^{-t} \left( \int_{-e^t}^{-1} |x|\nu(dx) +
\int_1^{e^t}|x|\nu(dx)\right)dt\\
&=\int_{-\infty}^{-1}\left(\int_{\log (-x)}^{\infty}e^{-t}dt\right)\nu(dx)
+ \int_1^{\infty}\left(\int_{\log x}^{\infty}e^{-t}dt\right)\nu(dx)\\
& = \int_{\{|x|>1\}} \frac{1}{|x|} \nu(dx) \leq \int_{\{|x|>1\}}\nu(dx) <\infty.
\end{align*}
Therefore,  we finally have
\begin{align*}
\lim_{A\to\infty}\int_{0}^{A}
\mC_{\law(Z_{1})}(e^{-t}z)dt, =
\int_{0}^{\infty}
\mC_{\law(Z_{1})}(e^{-t}z)dt
=\int_{0}^{1}
\mC_{\law(Z_{1})}(wz)\frac{dw}{w}, 
\quad z\in\C^{-}.
\end{align*}

We next want to show
\begin{eqnarray}\label{change}
\frac{d}{dz} \int_{0}^{1} \mathcal{C}_{\mathcal{L}(Z_{1})} (wz) \frac{dw}{w} 
=\int_{0}^{1} \frac{d}{dz} \mathcal{C}_{\mathcal{L}(Z_{1})} (wz)\frac{dw}{w} .
\end{eqnarray}
To do it, we start with
\begin{align}\label{derivative-free-cum}
\frac{d}{dz}\mathcal{C}_{\mathcal{L}(Z_{1})} (wz) 
= \eta w + 2aw^{2}z + \frac{d}{dz}\left(\int_{\R}\left( \frac1{1-xwz} - 1 - xwz{\bf 1}_{[-1,1]}(x)\right) \nu(dx)\right).
\end{align}
The first two terms have no problems. We thus omit them.
As to the third term, we want to show
\begin{align}
&\int_0^1\frac{d}{dz}\left(\int_{\R}\left( \frac1{1-xwz} - 1 - xwz{\bf 1}_{[-1,1]}(x)\right) \nu(dx)
\right)
\frac{dw}{w}\nonumber\\
 &\q\q = \int_0^1\int_{\R}\frac{d}{dz}\left( \frac1{1-xwz} - 1 - xwz{\bf 1}_{[-1,1]}(x)\right) \nu(dx)\frac{dw}{w}.
\label{change2}
\end{align}
We have 
\begin{eqnarray}\label{deriv}
\frac{d}{dz}\left( \frac1{1-xwz} - 1 - xwz{\bf 1}_{[-1,1]}(x)\right)
= \frac{xw}{(1-xwz)^{2}} - xw\mathbf{1}_{[-1,1]}(x).
\end{eqnarray}
Let $K$ be an arbitrary compact subset in $\C^{-}$.
We shall show that $\left|\frac{xw}{(1-xwz)^{2}} - xw\mathbf{1}_{[-1,1]}(x)\right|$ is dominated by an integrable function, which is locally independent of $z$.
On $\{(w,x)|w\in[0,1], x\in[-1,1]\}$, we have
\begin{align*}
&\left|\frac{xw}{(1-xwz)^{2}} - xw\mathbf{1}_{[-1,1]}(x)\right| \leq
\frac{\left|xw(1-(1-xwz)^{2})\right|}{\left|1-xwz\right|^{2}}\\
&\leq \frac{2 w^{2} |z||x|^{2} +w^{3} |z|^{2} |x|^{3}  }{d_{1}^{2}}
\leq \frac{ (2|z| + |z|^{2})}{d_{1}^{2}}w|x|^{2}
\leq D_{K}w|x|^{2}, \q \forall z \in K,
\end{align*}
where $d_{1}:=\min_{0<t \le 1, z\in K}|1-tz|^2$ and $D_{K}:=(\max_{z\in K}(2|z| + |z|^{2}))/d_1$.
On $\{(w,x)|w\in[0,1], x\in\R\backslash[-1,1]\}$, we have
\begin{align*}
\left|\frac{xw}{(1-xwz)^{2}} - xw\mathbf{1}_{[-1,1]}(x)\right| 
= \left|\frac{xw}{(1-xwz)^{2}}\right| \leq \frac{w}{|x||x^{-1} - wz|^{2}} \leq \frac{w}{d_{2}^{2}}, \q \forall z \in K,
\end{align*}
where $d_{2}^{2} = \min_{z\in K, x \in \R\backslash[-1,1]}\{ |x| |x^{-1}-wz|^{2}\}$.
Therefore, for any compact subset $K$ in $\C^{-}$, there exists a positive constant $C_3>0$ such that, 
\begin{align*}
\left|\frac{xw}{(1-xwz)^{2}} - xw\mathbf{1}_{[-1,1]}(x)\right|\leq C_3 w(|x|^{2}\wedge 1),\quad  z \in K, w\in[0,1], x\in\R,
\end{align*}
where the function on the right hand side is a dominating function that
is integrable by $\nu(dx)\frac{dw}{w}$ on $\{(w,x)|w\in[0,1], x\in\R\}$.
Thus,  by the theorem on the
differentiation under the integral sign, we can finally conclude that
\begin{align*}
&\quad\,\, z \frac{d}{dz} \mathcal{C}_{\mathcal{L}(X)} (z) 
=z \frac{d}{dz} \int_{0}^{1} \mathcal{C}_{\mathcal{L}(Z_{1})} (wz) \frac{dw}{w} \\
&\q\q =z \int_{0}^{1} \frac{d}{dz} \mathcal{C}_{\mathcal{L}(Z_{1})} (wz)\frac{dw}{w} 
=z \int_{0}^{1} \frac{d\mathcal{C}_{\mathcal{L}(Z_{1})}}{dz} (wz)dw \\
&\q\q =\int_{0}^{1} \frac{d}{dw}\mathcal{C}_{\mathcal{L}(Z_{1})}(wz)dw 
=\mathcal{C}_{\mathcal{L}(Z_{1})}(z),
\end{align*}
completing the proof.
\qed

\vskip 3mm

\begin{thm}\label{final}
If the free characteristic triplet of $\law(X)$ is $(a_{X},\nu_{X},\eta_{X})$, then
\begin{align*}
\mC_{\law(Z_{1})}(z) 
= \eta_{X} z + 2a_{X}z^{2} 
+  z\int_{\R} \left(\frac{x}{(1-zx)^{2}} - x\mathbf{1}_{[-1,1]}(x)\right)\nu_{X}(dx),\quad z \in \C^{-}.
\end{align*}
\end{thm}

{\it Proof.}
Note that \eqref{derivative-free-cum} with $w=1$ is
\begin{align*}
\frac{d}{dz}\mathcal{C}_{\mathcal{L}(X)} (z) 
= \eta_{X} + 2a_{X}z + \frac{d}{dz}\left(\int_{\R}\left( \frac1{1-xz} - 1 - xz{\bf 1}_{[-1,1]}(x)\right) \nu_{X}(dx)\right) ,
\q z\in\C^-.
\end{align*}
The validity of the interchange of derivative and integral above has already been proved
in the proof of Theorem~\ref{thmformula}.
Thus, \eqref{formula2} and \eqref{deriv} with $w=1$ imply the statement of this theorem.
\qed


\vskip 10mm
\section{Examples}
In this section, we shall discuss Theorem \ref{thmformula} or \ref{final} with explicit examples.
\begin{ex}[Semicircle distributions]
If the free characteristic triplet of $\mu\in ID(\bp)$ is $(a, 0, \eta)$,
$\mu$ is called the semicircle distribution $\mathbf{w}(\eta ,a)$, and it is freely
selfdecomposable.  (See \cite{BT06}.)
As we have seen in Theorem~\ref{thmformula}, the free cumulant of 
the background deriving free L\'evy process $\{Z_{t}\}$ of $\mathbf{w}(\eta, a)$ is
$$
\mC_{\law(Z_{1})} (z)= \eta z + 2az^{2}, \q z \in \C^{-}.
$$
Thus, $\{Z_{t}\}$ is the free Brownian motion which is a free L{\'e}vy process with the marginal distribution $\mathbf{w}(\eta ,2a)$ at time $t=1$.
\end{ex}

\begin{ex}[Free gamma process]
Let $t\geq 0$ and $c>0$.
Assume that $\gm(t,c) \in \mathcal{P}(\R)$ has the following free cumulant transform:
\begin{align*}
\mC_{\gm(t,c)} (z)
= \int_{(0,4c)} \left(\frac{1}{1-xz}-1\right)\frac{t\sqrt{x(4c-x)}}{2\pi x^{2}} d x
= \frac{t(1-\sqrt{1-4c z})}{2},\quad z \in \C^{-}.
\end{align*}
From this free cumulant transform we can obtain its Cauchy transform and probability density function 
through the definitions of the free cumulant transform and the Stieltjes inversion formula as follows, although we omit the detailed calculations here:
$$
G_{\gm(t,c)}(z) = \frac{z(t+2)-c - \sqrt{z^{2}-2cz(t+2) + c^{2}}}{2z^{2}}
$$
and
$$
f_{\gm(t,c)}(x) = \frac{t}{2\pi x^{2}}
\sqrt{\left(x - \alpha_{-}\right)\left(\alpha_{+} -x\right)}
\mathbf{1}_{[\alpha_{-},\alpha_{+}]}(x), 
$$
where $\alpha_{\pm} = c((t+2) \pm 2\sqrt{t+1})$.

$\gm(t,c)$ is called the (non-centered) free gamma distribution with shape parameter $t$ and scale parameter $c$ 
from the view of free Meixner, family, (see \cite{SY01, A03, BB06}).
Since the measure 
$\nu_{t,c}(dx) = \frac{t\sqrt{x(4c-x)}}{2\pi x^{2}}\mathbf{1}_{(0,4c)}(x)dx$ satisfies the condition of the L{\'e}vy measure for $\gm(t,c)$ and 
$l(x)=\frac{t\sqrt{x(4c-x)}}{2\pi x}\mathbf{1}_{(0,4c)}(x)$ is decreasing on $(0,\infty)$ and vanishes on $(-\infty,0)$, $\gm(t,c)$ is freely selfdecomposable. 
Note that it is not an image of the classical gamma distributions by the Bercovici-Pata bijection, and if $t=1$ and $c=1$, then the function $l(x)$ is the density function of the standard Marchenko-Pastur distribution (that is, the standard free Poisson distribution). 
The L{\'e}vy measure of the classical gamma distribution with parameter $(t,c)$ is of the form $\frac{t e^{-cx}}{x}\mathbf{1}_{(0,\infty)}(x)dx$.
Note that if $t=1$ and $c=1$ then the function $e^{-x}\mathbf{1}_{(0,\infty)}(x)$ is the density function of the standard classical exponential distribution.

In classical probability theory, as it is pointed out in Example 16.3 (iv) of \cite{S99}, 
an $H$-selfsimilar additive process
based on the exponential distributions has exponential distributions for all 
$t\geq 0$.
From this point of view, we can introduce free exponential distributions. 
We consider the laws $\{\mu_{t}, t\ge 0\}$ of an $H$-selfsimilar free additive process $\xt$ such that
$\law(X_1)=\gm(1,1)=:\mu$ has the following free cumulant transform:
\begin{align*}
\mC_{\mu_{t}} (z) 
= \mC_{\mu} (t^{H} z)
= \frac{1-\sqrt{1-4 t^{H} z}}{2}
\end{align*}
and 
the laws $\{\rho_{t}, t\ge 0\}$ of the free  L\'evy process $\yt$ such that
$\law(Y_1)=\mu$ has the following free cumulant transform:
\begin{align*}
\mC_{\rho_{t}} (z) 
= t \mC_{\mu} (z)
= \frac{t(1-\sqrt{1-4z})}{2}.
\end{align*}
As in classical probability theory,
we should refer to the free L\'evy process $\yt$ as the standard free gamma process.

Finally, we want to obtain the background driving free L\'evy process of the free selfdecomposable random variable $X_{1}$.
 By Theorem \ref{thmformula}, 
\begin{align*}
\mC_{\law(Z_{1})} (z) 
&=-\frac{-z}{\sqrt{1 - 4z}}.
\end{align*}
Applying the formula in p.304 of \cite{SSV}
\begin{align*}
\frac{\lambda}{\sqrt{\lambda + a}} 
= \int_{(a,\infty)} \frac{\lambda}{\lambda + u}\cdot\frac{1}{\pi\sqrt{u-a}} du
\end{align*}
to the above free cumulant, we obtain
\begin{align*}
\frac{z}{\sqrt{1-4z}} 
&= - \frac{1}{2} \frac{-z}{\sqrt{-z+1/4} }
= - \frac{1}{2} \int_{(1/4,\infty)} \frac{-z}{-z + u}\cdot\frac{1}{\pi\sqrt{u-1/4}}du\\
&= \frac{1}{2} \int_{(1/4,\infty)} \frac{u}{-z + u}-1\frac{du}{\pi\sqrt{u-1/4}}\\
&= \frac{1}{2} \int_{(0,4)} \left(\frac{1}{1-xz}-1\right) \frac{2}{\pi x \sqrt{x(4-x)}}dx.
\end{align*}
Therefore the free L\'evy measure of $Z_{1}$ is $\frac{dx}{\pi x \sqrt{x(4-x)}}$.
\end{ex}

\begin{rem}
Other free gamma distributions, which are the images of the classical gamma distributions under the Bercovici-Pata bijection, were studied by Haagerup and Thorbj{\o}rnsen \cite{HT14}. It follows from their L{\'e}vy measures that they are also freely selfdecomposable, and thus we can construct  $H$-selfsimilar free additive processes based on them. But we cannot get explicit forms of their densities, although asymptotic properties can be obtained from the results in \cite{HT14}.
\end{rem}

\begin{ex}
We consider $\mu(p)\in ID(\boxplus), -1<p<1$, whose free cumulant transform is as follows:
\begin{align*}
\mathcal{C}_{\mu(p)}(z) 
&= 1- (1-z)^{p}\\
&= \int_{\R\backslash\{0\}}\left(\frac{1}{1-xz}-1\right) \left(\frac{\sin (p\pi)}{\pi x}\right)\left(\frac{1-x}{x}\right)^{p} 
\mathbf{1}_{(0,1)}(x)d x.
\end{align*}
Especially, when $0<p<1$, it is freely selfdecomposable. 
Thus we can construct $\{X_{t}\}$ and $\{Z_{t}\}$.
The free cumulant transform of $Z_{1}$ is 
\begin{align*}
\mathcal{C}_{\law(Z_{1})}(z) 
= pz(1-z)^{p-1}.
\end{align*}
Applying the formula in p.305 of \cite{SSV} which is
\begin{align*}
\frac{\lambda}{(\lambda + a)^{\alpha}} 
= \lambda\sin(\alpha \pi)\int_{(a,\infty)} \frac{1}{\pi(\lambda + u)(u-a)^{\alpha}}du, \text{ for } a>0, \alpha\in(0,1),
\end{align*}
we obtain 
\begin{align*}
\mathcal{C}_{\law(Z_{1})}(z) = pz(1-z)^{p-1}
= p \sin((1-p)\pi)\int_{(0,1)} \left(\frac{1}{1-xz}-1\right) \frac{ 1}{\pi x^{1+p}(1-x)^{1-p}}dx.
\end{align*}
\end{ex}

\begin{ex}[Fuss Catalan distribution $\mu(p,p)$] 
We consider, for $1<p<2$, a probability measure $\mu(p,p)$, whose free cumulant transform is as follows:
\begin{align*}
\mathcal{C}&_{\mu(p,p)}(z) 
= pz + (z + 1)^{p} -1 \\
&= pz + \int_{\R\backslash\{0\}}\left(\frac{1}{1-xz}-1-zx\mathbf{1}_{[-1,1]}(x)\right) \left(-\frac{\sin (p\pi)}{\pi |x|}\right)\left(\frac{1+x}{-x}\right)^{p} 
\mathbf{1}_{(-1,0)}(x)d x.
\end{align*}
It is called the Fuss-Catalan distribution with parameter $(p,p)$. See \cite{MSU2020}.
From its  L\'evy measure, $\mu(p,p)$ is in $L(\boxplus)$.
We can obtain the free cumulant transform of $Z_{1}$ from the similar computation of the above example as follows:
\begin{align*}
\mC&_{\law(Z_{1})} (z) =pz +pz(z+1)^{p-1} \\
& = \int_{\R\backslash\{0\}}\left(\frac{1}{1-xz}-1-zx\mathbf{1}_{[-1,1]}(x)\right) \left(-\frac{\sin (p\pi)}{\pi |x|}\right)(1+x)^{p-1}(-x)^{-p} 
\mathbf{1}_{(-1,0)}(x)d x.
\end{align*}
\end{ex}

\vskip 10mm

\section*{Acknowledgements}

The authors would like to thank a referee for his/her carefully reading the manuscript
and giving a lot of valuable comments, which led to an improvement of this paper.
Especially, his/her comment on the selfsimilarity is very much appreciated.
NS was partially supported by JSPS Kakenhi 19H01791, 19K03515, JPJSBP120209921, JPJSBP120203202.


\end{document}